\documentclass{IEEEtran}
\ifCLASSINFOpdf
\else
\fi
\usepackage[cmex10]{amsmath}
\usepackage{amsfonts}
\usepackage{xcolor}
\usepackage{rotating}
\usepackage{graphicx,multicol}
\usepackage{amssymb}
\usepackage{bm}
\usepackage{makecell}
\usepackage{hyperref}

\DeclareMathOperator{\real}{Re}
\DeclareMathOperator{\imag}{Im}

\DeclareMathOperator{\var}{var}
\DeclareMathOperator{\subjectto}{s.t.}

\begin{document}
% \title{Families of Second Order Cone Constraints for KVL in 3-Cycles in Optimal Power Flow
\title{Notes on BIM and BFM  Optimal Power Flow With Parallel Lines and Total Current Limits
\author{Frederik~Geth,~\IEEEmembership{Member,~IEEE,}
        Bin Liu,~\IEEEmembership{Member,~IEEE}
\thanks{F. Geth and B. Liu are with the Energy Systems program, CSIRO Energy, Newcastle NSW, Australia (e-mail: frederik.geth@csiro.au, brian.liu@csiro.au).}
\thanks{This paper was inspired by discussion on the issue tracker of PowerModels: https://github.com/lanl-ansi/PowerModels.jl/pull/286}}}

%% To specify the authors when (number of affiliations > 2)
% \author{\IEEEauthorblockN{Author n.1\IEEEauthorrefmark{1},
% Author n.2\IEEEauthorrefmark{2},
% Author n.3\IEEEauthorrefmark{3}, 
% Author n.4\IEEEauthorrefmark{3} and
% Author n.5\IEEEauthorrefmark{4}}
% \IEEEauthorblockA{\IEEEauthorrefmark{1} Department Name of Organization A\\
% Name of the organization A,
% Address A\\ Emails if wanted}
% \IEEEauthorblockA{\IEEEauthorrefmark{2} Department Name of Organization B\\
% Name of the organization B,
% Address B\\ Emails if wanted}
% \IEEEauthorblockA{\IEEEauthorrefmark{3} Department Name of Organization C\\
% Name of the organization C,
% Address C\\ Emails if wanted}
% \IEEEauthorblockA{\IEEEauthorrefmark{4}Department Name of Organization D\\
% Name of the organization D,
% Address D\\ Emails if wanted}
% }

\maketitle

\newcommand{\eqdef}{\mathrel{\overset{\makebox[0pt]{\mbox{\normalfont\tiny\sffamily def}}}{=}}}

\begin{abstract}
The second-order cone relaxation of the branch flow model (BFM) and bus injection model (BIM) variants of optimal power flow are well-known to be equivalent for radial networks. 
In this work we show that in meshed networks with parallel lines, BIM dominates BFM, and propose novel constraints to make them equivalent in general.
Furthermore, we develop an improvement to the second-order cone relaxations of optimal power flow, adding novel and valid linear constraints on the lifted current expressions. 
We develop two simple test cases to highlight the advantages of the proposed constraints.
These novel constraints tighten the second-order cone relaxation gap on test cases in the `PG Lib' optimal power flow benchmark library, albeit generally in limited fashion. 
% Throughout, we pay special attention to test cases in PG Lib where physical intuitions don't hold. For instance, from first principles, we expect branch series resistance, series reactance and shunt susceptance to be positive, as well as generator power output. Nevertheless, in PG Lib, there are numerous cases were such assumptions are false. 
\end{abstract}


\section{Introduction}

Convex relaxation is a powerful technique to enable global optimization of optimal power flow (OPF).
When combined with techniques such as optimality-based bound tightening  and cut generation, researchers have been able to prove global optimality of OPF benchmarks of nontrivial size \cite{Xu2021,Gopinath2020,Sundar2018a}.
Generally, one wants very tight (root-node) relaxations, that are minimally expensive to evaluate, and the second-order cone (SOC) relaxations have proven themselves offering exactly this trade-off.

In the context of OPF, there have been two competing frameworks to obtain SOC relaxations: the branch flow model (BFM) and the bus injection model (BIM). 
The BIM develops expressions using the admittance form of Ohm's law, allowing all current variables to be eliminated in favour of voltage variables; the BFM still retains a variable for the series current to solve the impedance form of Ohm's law. 
The BIM and BFM  SOC relaxations were shown to be equivalent for radial networks by Low \cite{Low2014}. Equivalence  means the relaxations have the same feasible set, and that we know the bijection. 

\subsection{Different Styles of Branch Models}
Generally, authors focus on canonical BIM and BFM relaxations that do not include branch shunt admittance into the branch model (i.e. use a series impedance instead of a $\Pi$-section). 
We note that from the physics perspective, this just means composing a circuit from series and shunt impedances, and therefore does not introduce any error. 
Nevertheless, from the perspective of optimization problem feasibility, we need to pay special attention to the bound semantics: the power flow bound applies to the \emph{total} power flowing into the $\Pi$-section, not just to the power through the \emph{series} impedance. 

\begin{figure}[tbh]
  \centering
    \includegraphics[width=0.60\columnwidth]{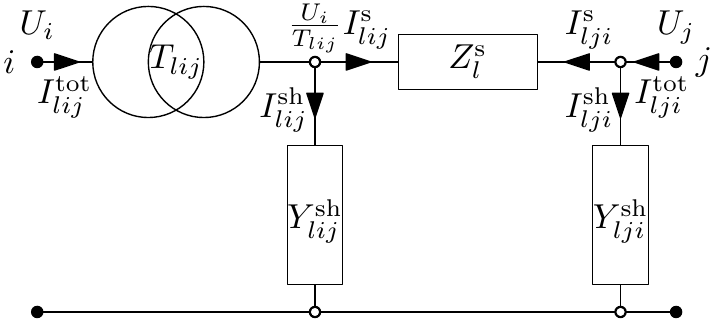}
  \caption{\textsc{Matpower}-style branch model: $\Pi$-section with ideal transformer}  \label{fig_linemodel_iv}
\end{figure}

Benchmark OPF libraries such as \textsc{PG Lib} \cite{Babaeinejadsarookolaee2019} use $\Pi$-model branches \emph{with an ideal transformer at the sending end} (Fig. \ref{fig_linemodel_iv}).
Coffrin et al. show  \cite{coffrin2015distflow} how to adapt the BIM and BFM SOC relaxations to incorporate shunts and the transformer ratio, and how to correctly apply the apparent power, voltage magnitude, and voltage angle difference bounds. 

Finally, we observe that,  \textsc{Matpower} \cite{Zimmerman2011}, and compatible OPF tools such as \textsc{PowerModels.jl} \cite{Coffrin2017}, by default do \emph{not} use branch \emph{current} limits, and neither do the benchmarks in \textsc{PG Lib}. 
Therefore, bounds on the lifted current variables in BFM are generally not included in the problem specification. 
Nevertheless, if buses have a nonzero lower voltage magnitude bound, it is possible to infer valid current bounds from the known apparent power limits.

\subsection{Issues With the Relaxations}\label{BIM_BFM}
\begin{figure}[tbh]
  \centering
    \includegraphics[scale=0.75]{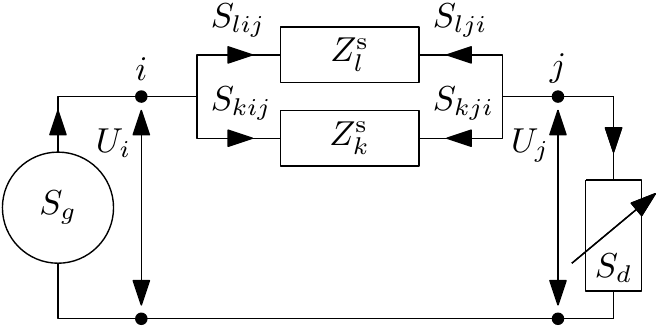}
  \caption{A 2-bus system with parallel lines.}  \label{fig-2-bus-system}
\end{figure}

We demonstrate next that the BIM SOC form is tighter than BFM when there are parallel lines. 
We therefore include an explicit branch index $l$, to disambiguate parallel lines (i.e. different lines between identical nodes $i$ and $j$).
Taking the simple 2-bus system in Fig.\,\ref{fig-2-bus-system} as an example, with branch $l$ and $k$ in parallel between buses $i$ and $j$, where branch impedances $Z_l$, $Z_k$ and the load set point $S_d=P_d+jQ_d$ are known parameters\footnote{For simplicity, operational constraints are not considered in this example.}. 
Voltages on bus $i$ and $j$ are $U_i, U_j$, the complex power flow through branch $l$ in the direction of $i$ to $j$ is $S_{lij}$.
We now solve for the generator dispatch value $S_g$.

The BIM SOC OPF is formulated as,
\begin{subequations}
\label{bim-2bus}
\begin{IEEEeqnarray}{LL}
\var &S_g, W_i, W_{j} ,W_{ij}, S_{lij}, S_{kij}, S_{lji}, S_{kji},\\
\min & f(\real(S_g)),\\
\subjectto& S_g=S_{lij}+S_{kij}, \quad S_{lji}+S_{kji}=-S_d,\\
&S_{lij} = Y_l^*( W_i - W_{ij}), \quad S_{kij} = Y_k^*( W_i - W_{ij}) \label{eq_bim_power_flow_from},\\
&S_{lji} = Y_l^*( W_j - W_{ij}^*), \quad S_{kji} = Y_k^*( W_j - W_{ij}^*),\\
&|W_{ij}|^2\le W_i W_j,
\end{IEEEeqnarray}
\end{subequations}
and the BFM  as,
\begin{subequations}
\label{bfm-2bus}
\begin{IEEEeqnarray}{LL}
\var & S_g, W_i, W_{j}, S_{lij}, S_{kij}, S_{lji}, S_{kji}, L_l , L_k, \\
\min & f(\real(S_g)),\\
\subjectto & S_g=S_{lij}+S_{kij}, \quad S_{lji}+S_{kji}=-S_d,\\
&S_{lij}+S_{lji}= Z_l L_l, \quad S_{kij}+S_{kji}=Z_k L_k ,\\
\label{bfmsoc2}
&|S_{lij}|^2\le W_i L_l,\quad  |S_{kij}|^2\le W_i L_k,\\
&W_j=W_i - Z^*_l S_{lij} - Z_l S^*_{lij} + |Z_l|^ 2L_l,\\
&W_j=W_i - Z^*_k S_{kij} - Z_k S^*_{kij} + |Z_k|^2 L_k.
\end{IEEEeqnarray}
\end{subequations}
It is noted that these sets have different numbers of independent variables, despite representing the same (nondegenerate) physical system. 
In the BIM, the bus-pair cross product variable $W_{ij}$ is shared\footnote{We note that the shared $W_{ij}$ variable approach is used to calculate the  SOC relaxation gap in PG Lib report \cite{Babaeinejadsarookolaee2019}.} amongst parallel branches between bus $i$ and $j$.
Conversely, in the BFM, there is an independent lifted current  variable  $L_l$ for each branch, regardless of whether it is in parallel or not.
Therefore, the feasible sets are not equivalent. 
An ill-informed way to make the BIM and BFM equivalent for parallel lines is to change the BIM formulation  with a branch-wise $W_{lij}$ instead of bus-pair-wise $W_{ij}$ variable:
\begin{IEEEeqnarray}{C}
\forall lij: W_{lij} \eqdef U_i U_j^* .
\end{IEEEeqnarray}
This obviously weakens the BIM relaxation for parallel lines.
Existing literature can be interpreted as taking this approach.

% Regarding a complex variable as two variables in real number domain, there will be 2 independent real variables in BIMTwo while the number is 4 for BFMTwo, implying their feasible regions are different, which, in the other words, means their optimal solutions may be different from each other. This is because the two lines use the common $W_{ij}$ in BIMTwo while two variables, i.e., $L_1$ and $L_2$, are introduced in BFMTwo. 

Alternatively,  we can strengthen the BFM for parallel lines.
We can re-write \eqref{eq_bim_power_flow_from} in terms of $W_{ij}$, 
\begin{IEEEeqnarray}{C}
W_{ij}=W_i-Z^*_l S_{lij}=W_i-Z^*_k S_{kij},
\end{IEEEeqnarray}
which implies
\begin{IEEEeqnarray}{C}
\label{bfmzs}
Z^*_lS_{lij}=Z^*_k S_{kij} \label{eq_parallel_lines_simple},
\end{IEEEeqnarray}
is the missing constraint that makes BIM and BFM equivalent again with parallel lines.
% , if \eqref{bfmzs} is added to BFMTwo, the two formulations will be equivalent to each other, and one of the SOC constraints in \eqref{bfmsoc2} can be removed.

% However, either adding bounds on $L_l$ or adding \eqref{bfmzs}for $\Pi$-section with transformer are much more complicated, which will be discussed in detail in the next section.   
We define apparent power branch flow limits on both sides,
\begin{IEEEeqnarray}{C}
|S_{lij}| \leq S_{l}^{\text{max}}, ~
|S_{lji}| \leq S_{l}^{\text{max}}.
\end{IEEEeqnarray}
We can now derive bounds on the left-hand side  and right-hand side (RHS) of \eqref{bfmsoc2} separately,
\begin{IEEEeqnarray}{C}
2 S_{l}^{\text{max}} \geq |S_{lij} + S_{lji}| =
|Z_l L_l| \leq Z_l \left(\frac{S_{l}^{\text{max}}}{U_i^{\text{min}}} \right)^2.
\end{IEEEeqnarray}

The RHS can be significantly tighter in certain situations, for instance ($Z_l, S_{l}^{\text{max}}, U_i^{\text{min}}$) = (0.1, 1, 0.9), where the ratio is 162:1. 
This suggests that implied current bounds can be powerful, and deserve further study.

\subsection{Scope and Contributions}
% While BIM and BFM SOC relaxations for equivalent for radial networks, the converse will be shown to  not be true: in meshed networks the extended BFM and BIM SOC relaxations \cite{coffrin2015distflow} do \emph{not} describe the same feasible set, in the context of networks with parallel lines. 
% We will develop examples where the canonical SOC BIM dominates the BFM \emph{and the other way around}, based on 2 and 3-bus test cases. 
% We then propose two improvements to restore the BIM and BFM SOC  equivalence, which improve the SOC relaxation gap for meshed networks in general.
In this work, we generalize the parallel line consistency constraint \eqref{eq_parallel_lines_simple} for the SOC BFM formulation subject to \textsc{Matpower}-style branch models.
Moreover, we  propose linear total current limit expressions for both BIM and BFM, and finally study the relaxation gaps numerically.

% Note that  the active and reactive losses are easy to derive explicitly in this variable space, i.e. 
% \begin{IEEEeqnarray}{C}
% P_{lij} + P_{lji} =  R_{l} L_l , \quad
% Q_{lij} + Q_{lji} =  X_{l} L_l .
% \end{IEEEeqnarray}
% Therefore the loss bounds are, for $r_{l}>0, x_{l}>0$,
% \begin{IEEEeqnarray}{C}
% 0 \leq P_{lij} + P_{lji} \leq  R_{l} \cdot \min((I_{lij}^{\text{max}})^2, (I_{lji}^{\text{max}})^2), \\
% 0 \leq Q_{lij} + Q_{lji} \leq  X_{l} \cdot \min((I_{lij}^{\text{max}})^2, (I_{lji}^{\text{max}})^2) ,
% \end{IEEEeqnarray}
% In nondimensionalized (per unit) quantities $|U_i| \approx 1$ and therefore, $|S_{lij}| \approx |I_{lij}|$. 
% Furthermore, in typical PG lib test cases,  $r_{l}<1$ pu, $x_{l}<1 $ pu, therefore, we can see that these bounds are often significantly tighter than \eqref{eq_loss_bounds_smax}. 

\section{BIM \& BFM for $\Pi$-sections with Transformer}
We now derive the BIM and BFM for \textsc{Matpower}-style $\Pi$-model branches (Fig.~\ref{fig_linemodel_iv}).
Note that this includes an idealized transformer, with complex tap ratio $T_{lij}$, at the sending end of the branch. 
The power and current bounds now apply at the terminals, not just to the series impedance. 
Therefore if the shunts are nonzero, the bound semantics are different from the canonical case.
We slightly generalize the \textsc{Matpower} branch model to support asymmetric shunts, and conductive shunts.
This allows for parameterizing a $\Gamma$-section as an edge case of a $\Pi$-section.

\subsection{Preliminaries}
The shunt at the from-side is  $Y^{\text{sh}}_{lij}$, at the to-side $Y^{\text{sh}}_{lji}$;
the series impedance is  $Z^{\text{s}}_{l}$;
the transformer ratio is $T_{lij}$.
The current through the series impedance is $I^{\text{s}}_{lij}$ in the direction of $i$ to $j$ and $I^{\text{s}}_{lji}$ for $j$ to $i$.
The current through the shunt at the from side is  $I^{\text{sh}}_{lij}$ and at the to side is  $I^{\text{sh}}_{lji}$.
Therefore the current divides over the series and shunt elements as,
\begin{IEEEeqnarray}{C}
 I^{\text{tot}}_{lij} =( I^{\text{s}}_{lij} + I_{lij}^{\text{sh}} )/T^*_{lij} , \quad
I^{\text{tot}}_{lji} = I^{\text{s}}_{lji} + I^{\text{sh}}_{lji}, \label{eq_kcl_branch} \\
I^{\text{s}}_{lij}  +  I^{\text{s}}_{lji}  = 0 . \label{eq_series_current_kcl}
\end{IEEEeqnarray}
The shunt currents $I_{lij}^{\text{sh}}, I_{lji}^{\text{sh}}$ depend on the voltage $U_i, U_j$,
\begin{IEEEeqnarray}{C}
I_{lij}^{\text{sh}} = Y^{\text{sh}}_{lij} \frac{U_i}{T_{lij}}, \quad I_{lji}^{\text{sh}} = Y^{\text{sh}}_{lji} U_j  . \label{eq_shunt_adm_iv}
\end{IEEEeqnarray}
Ohm's law  between buses $i$ and $j$ through branch $l$ is,
\begin{IEEEeqnarray}{C}
U_j = \frac{U_i}{T_{lij}} - Z^{\text{s}}_l I_{lij}^{\text{s}}. \label{eq_ohms_iv}
\end{IEEEeqnarray}
The (total) power flow variables are,
\begin{IEEEeqnarray}{C}
S_{lij}^{\text{tot}} \eqdef U_i (I_{lij}^{\text{tot}})^*, \quad
S_{lji}^{\text{tot}} \eqdef U_j (I_{lji}^{\text{tot}})^*, \label{eq_complex_power_total_def} \label{eq_tot_power_def}
\end{IEEEeqnarray}
and the apparent power limits apply,
\begin{IEEEeqnarray}{C}
|S_{lij}^{\text{tot}} |^2\leq (S_{l}^{\text{max}})^2,
|S_{lji}^{\text{tot}} |^2 \leq (S_{l}^{\text{max}})^2. \label{eq_tot_power_bounds}
\end{IEEEeqnarray}
We define the power flow through the series impedance, 
\begin{IEEEeqnarray}{C}
S_{lij}^{\text{s}} \eqdef \frac{U_i}{T_{lij}} (I_{lij}^{\text{s}})^* \label{eq_series_power_from_def}, \quad
S_{lji}^{\text{s}} \eqdef U_j (I_{lji}^{\text{s}})^* . \label{eq_series_power_to_def} 
\end{IEEEeqnarray}
Note that the sending end series power flow variable uses ${U_i}/{T_{lij}}$ -- not $U_i$ -- as the voltage variable.
We can substitute $ S_{lij}^{\text{tot}}$ for $S_{lij}^{\text{s}} $, by substituting \eqref{eq_kcl_branch} into \eqref{eq_tot_power_def},
\begin{IEEEeqnarray}{C}
 S_{lij}^{\text{tot}} = (Y^{\text{sh}}_{lij})^* \frac{W_i}{|T_{lij}|^2} + S_{lij}^{\text{s}}, ~
 S_{lji}^{\text{tot}} = (Y^{\text{sh}}_{lji})^* W_j + S_{lji}^{\text{s}} \label{eq_total_series_power_link}.
\end{IEEEeqnarray}
The lifted series current variable $L_l^{\text{s}}$ is,
\begin{IEEEeqnarray}{C}
0 \leq L_l^{\text{s}} \eqdef I_{lij}^{\text{s}} ( I_{lij}^{\text{s}} )^* = I_{lji}^{\text{s}} ( I_{lji}^{\text{s}} )^*. \label{eq_lifted_current_var_def}
\end{IEEEeqnarray}
The lifted voltage variable $W_i$ is,
\begin{IEEEeqnarray}{C}
(U_i^{\text{min}})^2 \leq W_i \eqdef U_i U_i^* \leq (U_i^{\text{max}})^2. \label{eq_wi_var_bounds}
\end{IEEEeqnarray}
The lifted voltage cross-product variable $W_{ij}$ is,
\begin{IEEEeqnarray}{C}
-U_i^{\text{max}} U_j^{\text{max}} \leq W_{ij} \eqdef U_i U_j^* \leq U_i^{\text{max}} U_j^{\text{max}}. \label{eq_wij_var_bounds}
\end{IEEEeqnarray}
The angle differences between adjacent buses are constrained through $\theta_{ij}^{\text{max}}$,
\begin{IEEEeqnarray}{C}
\label{anglelimit}
- \pi/4 \leq \theta_{ij}^{\text{min}} \leq \angle U_i - \angle U_j =  \angle U_i U_j^*\leq \theta_{ij}^{\text{max}} \leq \pi/4. \phantom{-}
\end{IEEEeqnarray}
In terms of $W_{ij}$, this becomes,
\begin{IEEEeqnarray}{C}
\tan (\theta_{ij}^{\text{min}} ) \real(W_{ij}) \leq \imag(W_{ij}) \leq \tan (\theta_{ij}^{\text{max}} ) \real(W_{ij}). \label{eq_wij_angle_diff_bounds}
\end{IEEEeqnarray}

\subsection{BFM SOC Formulation}
Ohm's law is obtained by multiplying \eqref{eq_ohms_iv} with its own conjugate, and substituting in the lifted variables,
\begin{IEEEeqnarray}{C}
% S_{lij}^{\text{tot}} + S_{lji}^{\text{tot}} = (y^{\text{sh}}_{lij})^* \frac{W_i}{|T_{lij}|^2}+ z^{\text{s}}_{l} L^{\text{s}}_l + (y^{\text{sh}}_{lji})^* W_j\\
W_j = \frac{W_i}{|T_{lij}|^2} - \left(Z_l^*  S^{\text{s}}_{lij} + \left(Z_l^*   S^{\text{s}}_{lij} \right)^* \right) + |Z_l|^2 L^{\text{s}}_l.
\end{IEEEeqnarray}
The SOC constraints are obtained by multiplying \eqref{eq_series_power_from_def} with its conjugate, and then relaxing to an inequality,
\begin{IEEEeqnarray}{C}
% S_{lij}^{\text{s}} = \frac{U_i}{T_{lij}} I^{\text{s}}_{lij} \\
|S_{lij}^{\text{s}}|^2 \leq \frac{W_i}{|T_{lij}|^2} L_l^{\text{s}} \label{eq_bfm_soc_const}
% , \quad
% |S_{lji}^{\text{s}}|^2 \leq {W_j} L_l^{\text{s}}
.
\end{IEEEeqnarray}
We obtain the loss balance by multiplying \eqref{eq_ohms_iv} with $(I^{\text{s}}_{lij})^*$
\begin{IEEEeqnarray}{C}
\label{bfmpowerloss}
S_{lij}^{\text{s}} + S_{lji}^{\text{s}} = Z^{\text{s}}_{l} L^{\text{s}}_l \label{eq_bfm_loss_series} .
\end{IEEEeqnarray}
% In the reals this becomes
% \begin{IEEEeqnarray}{C}
% (T^{\text{re}}_l P_{lij}^{\text{s}} + T^{\text{im}}_lQ_{lij}^{\text{s}}) /|T_{lij}|^2 + P_{lji}^{\text{s}} = R^{\text{s}}_{l} L^{\text{s}}_l  \\
% (T^{\text{re}}_l Q_{lij}^{\text{s}} - T^{\text{im}}_l P_{lij}^{\text{s}}) /|T_{lij}|^2  + Q_{lji}^{\text{s}} = X^{\text{s}}_{l} L^{\text{s}}_l 
% \end{IEEEeqnarray}

\subsection{BIM SOC Formulation}
We rewrite Ohm's law \eqref{eq_ohms_iv} in admittance form
\begin{IEEEeqnarray}{C}
I^{\text{s}}_{lij} =  Y^{\text{s}}_{l} \left(\frac{U_i}{T_{lij}} - U_j \right), \quad
I^{\text{s}}_{lji} =  Y^{\text{s}}_{l} \left( U_j - \frac{U_i}{T_{lij}} \right)\label{eq_ohms_adm_iv}.
\end{IEEEeqnarray}
% \begin{IEEEeqnarray}{C}
% (I^{\text{s}}_{lij})^* =  (Y^{\text{s}}_{l})^* (\frac{U_i^*}{T_{lij}^*} - U_j^*) \\
% (I^{\text{s}}_{lji})^* =  (Y^{\text{s}}_{l})^* ( U_j^* - \frac{U_i^*}{T_{lij}^*}) 
% \end{IEEEeqnarray}
We take the conjugate of \eqref{eq_ohms_adm_iv} and multiply with $U_i/T_{lij}, U_j$ to obtain,
\begin{IEEEeqnarray}{C}
S_{lij}^{\text{s}}  = (Y^{\text{s}}_l)^* \left( \frac{W_i}{|T_{lij}|^2} -  \frac{W_{ij}}{T_{lij}} \right) \label{eq_power_from_bim},\\
S_{lji}^{\text{s}}  = (Y^{\text{s}}_l)^* \left(W_j - \frac{W_{ji}}{T_{lij}^*} \right) \label{eq_power_to_bim}.
\end{IEEEeqnarray}
We get the total flow expressions by substituting the above into \eqref{eq_total_series_power_link}.
The SOC constraint, to link everything together, is,
\begin{IEEEeqnarray}{C}
|W_{ij}|^2 \leq W_i W_j \label{eq_wij_soc_relax}.
\end{IEEEeqnarray}

\subsection{Valid Current Bounds}

Given the total power limit $S_{l}^{\text{max}}$, we can derive a valid bound on the series current, at the sending end,
\begin{IEEEeqnarray}{C}
 |I_{lij}^{\text{tot}}| \leq S_{l}^{\text{max}}/ U_{i}^{\text{min}}, \quad
 |I_{lij}^{\text{sh}}| \leq |Y^{\text{sh}}_{lij}| U_{i}^{\text{max}} /|T_{lij}|,\\
 \implies |I_{lij}^{\text{s}}| \leq I_{lij}^{\text{s,max}} \eqdef |T_{lij}| |I_{lij}^{\text{tot}}| + |I_{lij}^{\text{sh}}| .
 \end{IEEEeqnarray}
and at the receiving end,
\begin{IEEEeqnarray}{C} \label{series-current-limit}
  |I_{lji}^{\text{tot}}| \leq S_{l}^{\text{max}}/ U_{j}^{\text{min}},
 |I_{lji}^{\text{sh}}| \leq |Y^{\text{sh}}_{lji}| U_{j}^{\text{max}},\\
 \implies  |I_{lji}^{\text{s}}| \leq I_{lji}^{\text{s,max}} \eqdef |I_{lji}^{\text{tot}}| + |I_{lji}^{\text{sh}}|. 
\end{IEEEeqnarray}
Finally, using \eqref{eq_series_current_kcl} we know $|I_{lij}^{\text{s}}| = |I_{lji}^{\text{s}}|$, and therefore
\begin{subequations}
\label{eq_series_current_bound}
\begin{IEEEeqnarray}{C}
\label{Ltot}
I_{l}^{\text{s,max}} \eqdef \min(I_{lij}^{\text{s,max}},I_{lji}^{\text{s,max}}), \implies
L^s_l \leq  (I_{l}^{\text{s,max}})^2. \label{eq_lifted_series_current_bound}
\end{IEEEeqnarray}
\end{subequations}
% We can now impose bounds on the total receiving/sending current, through,
% \begin{IEEEeqnarray}{C} \label{eq-total-Slimit}
% |T_{l}|^2 L^{\text{tot}}_{lij} = L^{\text{s}}_{l} + Y^{\text{sh}}_{lij} S^{\text{s}}_{lij} + (Y^{\text{sh}}_{lij}S^{\text{s}}_{lij})^* + |Y^{\text{sh}}_{lij}|^2 \frac{W_i}{|T_{lij}|^2}  \label{eq_tot_current_var_lifted} 
% \end{IEEEeqnarray}
% It is therefore possible to put down a linear constraint on the total current magnitudes as \eqref{Ltot} for the BFM SOC model.
We can  still enforce a total current limit in this variable space. We substitute \eqref{eq_shunt_adm_iv} into \eqref{eq_kcl_branch}, 
\begin{IEEEeqnarray}{C}\label{eq-total-current-01}
T_{lij}^*I^{\text{tot}}_{lij} = \left( I^{\text{s}}_{lij} + Y^{\text{sh}}_{lij} \frac{U_i}{T_{lij}} \right), 
I^{\text{tot}}_{lji} = I^{\text{s}}_{lji} + Y^{\text{sh}}_{lji} U_j.
\end{IEEEeqnarray}
Note that this procedure guarantees that  either of the total current limits is binding before the series current. 

Now we multiply this expression with its own conjugate and perform the variable substitutions, 
\begin{IEEEeqnarray}{C} \label{eq-total-Slimit}
|T_{lij}|^2 L^{\text{tot}}_{lij} = L^{\text{s}}_{l} + Y^{\text{sh}}_{lij} S^{\text{s}}_{lij} + (Y^{\text{sh}}_{lij}S^{\text{s}}_{lij})^*+ |Y^{\text{sh}}_{lij}|^2 \frac{W_i}{|T_{lij}|^2}\nonumber\\
\leq |T_{lij}|^2\left(\frac{S_{l}^{\text{max}}}{ U_{i}^{\text{min}}}\right)^2.
\label{eq_tot_current_var_lifted} 
\end{IEEEeqnarray}
Similarly, the receiving end total lifted current variable is,
\begin{IEEEeqnarray}{C}
L^{\text{tot}}_{lji} = L^{\text{s}}_{l} + Y^{\text{sh}}_{lji} S^{\text{s}}_{lji} + (Y^{\text{sh}}_{lji}S^{\text{s}}_{lji})^* + |Y^{\text{sh}}_{lji}|^2 {W_j}
\leq \left(\frac{S_{l}^{\text{max}}}{ U_{j}^{\text{min}}}\right)^2. \nonumber\\
\label{eq_tot_current_var_lifted_to}
\end{IEEEeqnarray}
For the BIM SOC model, introducing such constraints is  more complicated noting that $L^s_l$ is not defined in the  BIM variable space. Therefore, we introduce linking equalities.
% Alternatively, they can be imposed via the power loss constraints, which is based on \eqref{bfmpowerloss} are presented as \eqref{eq_loss_bound_1}-\eqref{eq_loss_bound_8} in the Appendix.
% We can now derive valid bounds for the series power flow, 
% \begin{IEEEeqnarray}{C}
% |S_{lij}^{\text{s}}| \leq   U_{i}^{\text{max}} I_{l}^{\text{s,max}}/|T_{lij}|~
% |S_{lij}^{\text{s}}| \leq    U_{j}^{\text{max}} I_{l}^{\text{s,max}}
% \end{IEEEeqnarray}

% Note that a valid bound on $L_l^{\text{s}}$ can now be imposed,
% \begin{IEEEeqnarray}{C} \label{series-current-limit-01}
% L_l^{\text{s}}  \leq (I_{l}^{\text{s,max}})^2. \label{eq_series_current_bound}
% \end{IEEEeqnarray}

% Similarly, the receiving end total lifted current variable is,
% \begin{IEEEeqnarray}{C}
% L^{\text{tot}}_{lji} = L^{\text{s}}_{l} + Y^{\text{sh}}_{lji} S^{\text{s}}_{lji} + (Y^{\text{sh}}_{lji}S^{\text{s}}_{lji})^* + |Y^{\text{sh}}_{lji}|^2 {W_j} \label{eq_tot_current_var_lifted_to}
% \end{IEEEeqnarray}

% It is therefore possible to put down a linear constraint on the total current magnitudes,
% \begin{IEEEeqnarray}{C}
% L^{\text{tot}}_{lij} \leq (I_{lij}^{\text{max}})^2, L^{\text{tot}}_{lji} \leq (I_{lji}^{\text{max}})^2
% \end{IEEEeqnarray}

% In terms of total power variables
% \begin{IEEEeqnarray}{C} \label{eq-total-Slimit_pq}
% |S_{lij}|^2 \leq W_i (I_{l}^{\text{s,max}})^2, \quad 
% |S_{lji}|^2 \leq W_j (I_{l}^{\text{s,max}})^2.
% \end{IEEEeqnarray}
% In terms of lifted series current variables

\subsection{Linking Equalities} \label{sec_linking_equalities}
Generally, the variable $L_l$ is only used in the BFM, and $W_{ij}$ only in  BIM. 
Nevertheless, as the relaxations are incomparable, as discussed in Section \ref{BIM_BFM}, it is useful to intersect them. 
Therefore we want to 1) derive the $W_{ij}$ variable as an expression of the natural BFM variables 2) derive the $L_l$ variable as an expression of the natural BIM variables. 

\subsubsection{BIM expressions in terms of $W_i, S_{lij}^{\text{s}}$ }
For the first case, we rewrite \eqref{eq_power_from_bim}-\eqref{eq_power_to_bim},
\begin{subequations}
\begin{IEEEeqnarray}{C}
W_{ij} = T_{lij} \frac{W_i}{|T_{lij}|^2}  - (Z^{\text{s}}_l)^* T_{lij} S_{lij}^{\text{s}}, \label{eq_wij_for_bfm_from}\\
W_{ji} = T_{lij}^* W_j  - (Z^{\text{s}}_l)^* T_{lij}^*. S_{lji}^{\text{s}}
\end{IEEEeqnarray}
\end{subequations}

We can now project the known bounds and the 4D cuts  \cite{Coffrin2015c} on $W_{ij}$ onto $S^{\text{s}}_{lij}, L^{\text{s}}_l$. 
The case of the parallel branches still needs to be considered.
We consider a case with two parallel branches $(l,i,j)$  and $(k,i,j)$, which must necessarily have identical  $W_{ij}$ values, and use  \eqref{eq_wij_for_bfm_from} to obtain,
\begin{IEEEeqnarray}{C}
\label{linking-01}
 T_{lij} \frac{W_i}{|T_{lij}|^2}  - (Z^{\text{s}}_l)^* T_{lij} S_{lij}^{\text{s}} =  T_{kij} \frac{W_i}{|T_{kij}|^2}  - (Z^{\text{s}}_k)^* T_{kij} S_{kij}^{\text{s}} .\nonumber \\ \label{eq_parallel_lines_bfm}
\end{IEEEeqnarray}
which generalizes \eqref{eq_parallel_lines_simple} to the \textsc{Matpower}-style branch model.
We need such a constraint to link the power flow variables through parallel lines pair-wise in the same way that $W_{ij}$ does.

% In the anti-parallel situation we have $(l,i,j)$ in parallel with $(k,j,i)$ and we use $W_{ij} = W_{ji}^*$
% \begin{IEEEeqnarray}{C}
% \footnotesize
% %  T_{lij} \frac{W_i}{|T_{lij}|^2}  - (Z^{\text{s}}_l)^* T_{lij} S_{lij}^{\text{s}} = 
% %  ( T_{kij}^* W_j  - (Z^{\text{s}}_k)^* T_{kij}^* S_{kji}^{\text{s}})^* \nonumber \\ 
%   T_{lij} \frac{W_i}{|T_{lij}|^2}  - (Z^{\text{s}}_l)^* T_{lij} S_{lij}^{\text{s}} = 
%   T_{kij} W_j^*  - Z^{\text{s}}_k T_{kij} (S_{kji}^{\text{s}})^* 
% \end{IEEEeqnarray}

\subsubsection{BFM expressions in terms of $W_i, W_{ij}$ }
For the second case, we substitute \eqref{eq_power_from_bim}-\eqref{eq_power_to_bim} into \eqref{eq_bfm_loss_series}, and divide by $Z_l$,
\begin{IEEEeqnarray}{C}\label{Llimit_01}
\footnotesize
\!\!L^{\text{s}}_l = |Y^{\text{s}}_l|^2 \!\! \left(\frac{W_i}{|T_{lij}|^2} - T_{lij}^* \frac{W_{ij}}{|T_{lij}|^2} - T_{lij} \frac{W_{ji}}{|T_{lij}|^2} + W_j \right)\label{eq_series_current_sq_bim}.
\end{IEEEeqnarray}
allowing us to place linear bounds on the total current magnitude when combined with \eqref{eq_tot_current_var_lifted}-\eqref{eq_tot_current_var_lifted_to}.

\section{Numerical Experiments}
We want to compare different formulations, in terms of gap w.r.t. the upper bound AC polar solution, as well as computation time, for:
\begin{enumerate}
    \item the canonical BIM relaxation used in \textsc{PG Lib} \cite{Babaeinejadsarookolaee2019} as implemented in \textsc{PowerModels.jl} \cite{Coffrin2017};

    \item the improved BFM relaxation with added implied total current limits  \eqref{eq_series_current_bound} and  parallel lines constraints  \eqref{eq_parallel_lines_bfm};
    
    \item the improved BIM relaxation with added the total current limits  \eqref{eq_series_current_bound}.
\end{enumerate}
We use \textsc{Ipopt} \cite{Wachter2006} as the SOC solver with \textsc{Mumps} as the linear solver. 

\subsection{Feasible Sets}

The KCL expression is the same for all variants,
\begin{IEEEeqnarray}{C}
\forall i: \sum_{lij \in \mathcal{T}} S^{\text{tot}}_{lij} + \sum_{d \in \mathcal{D} } S_d - \sum_{g \in \mathcal{G}} S_g = 0 \label{eq_kcl_power}.
\end{IEEEeqnarray}
and the generator output lies in a PQ space $\mathcal{S}_g$,
\begin{IEEEeqnarray}{C}
 S_g \in \mathcal{S}_g \label{eq_gen_freedom}.
\end{IEEEeqnarray}
We minimize the generation cost, using linear and quadratic coefficients $c_{g,1}, c_{g,2}$ w.r.t active power output $P_g = \real(S_g)$,
\begin{IEEEeqnarray}{C}
\min \sum_{g \in \mathcal{G}} \left( c_{g,1} P_g + c_{g,2} (P_g)^2 \right). \label{eq_objective}
\end{IEEEeqnarray}

\subsubsection{\textsc{Shared} constraints}
\begin{itemize}
\item generation cost objective \eqref{eq_objective};
\item complex power flow variables \eqref{eq_complex_power_total_def} or \eqref{eq_series_power_to_def}\footnote{The feasible sets are equivalent  either series or total complex power variables; nevertheless in the BFM the Ohm's law expressions are more simple in series power variables, which simplifies implementation.};
    \item apparent power flow bounds \eqref{eq_tot_power_bounds};
\item lifted voltage variable with bounds \eqref{eq_wi_var_bounds};
    \item KCL with constant power loads \eqref{eq_kcl_power};
    \item generator output \eqref{eq_gen_freedom};
    \item shunt power loss in branch \eqref{eq_total_series_power_link};
    \item lifted (4D) nonlinear cuts \cite{coffrin2015distflow}.
\end{itemize}

\subsubsection{\textsc{BIM Canonical}}
\textsc{Shared} + 
\begin{itemize}
\item lifted voltage crossproduct variable bounds \eqref{eq_wij_var_bounds};
    \item series power flow in branch \eqref{eq_power_from_bim}, \eqref{eq_power_to_bim};
    \item SOC linking constraint \eqref{eq_wij_soc_relax};
    \item voltage angle difference constraints \eqref{eq_wij_angle_diff_bounds}.
\end{itemize}

\subsubsection{\textsc{BFM Canonical}}
\textsc{Shared} +
\begin{itemize}
\item lifted series current variable \eqref{eq_lifted_current_var_def}, \eqref{eq_lifted_series_current_bound}
    \item series power flow in branch \eqref{eq_bfm_loss_series};
    \item SOC linking constraint \eqref{eq_bfm_soc_const};
    \item voltage angle difference constraints \eqref{eq_wij_angle_diff_bounds} + \eqref{eq_wij_for_bfm_from};
    \item linking variable for 4D nonlinear cuts \eqref{eq_wij_for_bfm_from}.
\end{itemize}

\subsubsection{\textsc{BIM Improved}}
\textsc{BIM Canonical} + 
\begin{itemize}
    \item total current bounds \eqref{eq_tot_current_var_lifted}, \eqref{eq_tot_current_var_lifted_to};
    \item lifted series current expression \eqref{eq_series_current_sq_bim}.
\end{itemize}

\subsubsection{\textsc{BFM Improved}}
\textsc{BFM Canonical} + 
\begin{itemize}
    \item total current bounds \eqref{eq_tot_current_var_lifted}, \eqref{eq_tot_current_var_lifted_to};
    \item parallel lines consistency \eqref{eq_parallel_lines_bfm}.
\end{itemize}

\subsection{Novel Small Test Cases}

We propose two simple test cases\footnote{m files available here: \url{https://doi.org/10.25919/znc4-5z12}} to highlight the improvements proposed in this article.
\subsubsection{case2\_parallel}
In this newly proposed 2-bus test case with 2 lines in parallel (Similar to Fig.\,\ref{fig-2-bus-system}, but with generator 2 added on bus 2,  parameters in Table\,\ref{tab_2_parallel_system}). 
The case is set up to have congestion on line 2. 
The generator on bus 2 is more expensive than on line 1, so preferably gen 1 would be dispatched. 
Without the parallel line constraint, it is possible to control the flow through the branches independently, and therefore gen 1 gets dispatched to supply all of the load, i.e. at a value of $P_{g=1}=1.146$. 
With the parallel line constraint, Kirchhoff's voltage law is correctly applied, and current is shared appropriately between the branches, and therefore both generators need to be dispatched $P_{g=1}=0.075$, $P_{g=2}=1.040$.

The BIM relaxation is exact (0 gap), matching the value of the ACOPF upper bound at 5.27360.  The BFM relaxation however has a gap of 78.264 \%, which goes to 0 when adding the parallel line constraints.

\begin{table}[h]
\scriptsize
 \begin{center}
   \caption{Parameters for the 2-bus parallel system, where the units for all parameters are p.u. except for $\Delta\theta_{12}^\text{min},\Delta\theta_{12}^\text{max}$ (in degree) and $c_1,c_0$, which are generator cost coefficients.}\label{tab_2_parallel_system}
   	\setlength{\tabcolsep}{2 pt}
	\renewcommand\arraystretch{1.2}
   \begin{tabular}{l l | l l}
     \hline
        Symbol      & Value & Symbol      & Value  \\
        \hline
        $Z_{l}, Y^{\text{sh}}_{lij}, Y^{\text{sh}}_{lji}$ & $0.065+j0.62, j0.225, j0.225$ &$S_{l}^\text{max}$ & $90$  \\
        $Z_{k},  Y^{\text{sh}}_{kij}, Y^{\text{sh}}_{kji}$ & $0.025-j0.75, j0.35, j0.35$ &$S_{k}^\text{max}$ & $0.5$  \\
        $g=2: c_{g,1}, c_{g,2}$ & 5.0, 0.0 & $\Delta\theta_{ij}^\text{min},\Delta\theta_{ij}^\text{max}$ & $-30^\circ,30^\circ$\\
        $g \in {1,2}: P_{g}^\text{min},P_{g}^\text{max}$ & $0.0,2.0$ &$U_i^\text{min}, U_i^\text{max}$ & $0.90,1.10$\\
        $g \in {1,2}: Q_{g}^\text{min},Q_{g}^\text{max}$ & $-10.0,1.0$  &$U_j^\text{min}, U_j^\text{max}$ & $0.90,1.10$\\
        $g=1: c_{g,1}, c_{g,2}$ & 1.0, 0.0  &  $S_{d}$ & $1.1+j0.4$\\
        
        \hline 
   \end{tabular}
 \end{center}
\end{table} 
\subsubsection{case2\_gap}
% In this test case\footnote{adapted from https://github.com/lanl-ansi/PowerModels.jl/blob/master/ \\
% test/data/matpower/case5\_gap.m}, the BIM relaxation  has a gap of 2.69 \%, but the proposed BFM relaxation has a gap of only 0.5915\%. 

\begin{figure}[tbh]
  \centering
    \includegraphics[scale=0.7]{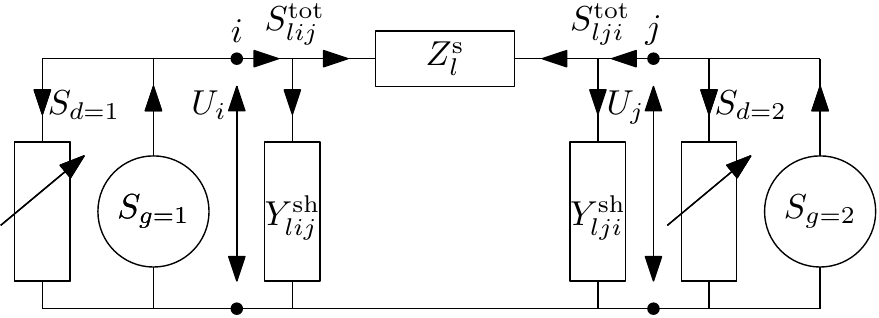}
  \caption{The 2-bus system to demonstrate the power of implied current limits.}  \label{fig_2_bus_system}
\end{figure}

% \subsection{Case2\_gap}
We fine tune a 2-bus 1-branch system to create a significant gap, to demonstrate the effectiveness of the implied constraints on the total current.
The circuit is presented in Fig.~\ref{fig_2_bus_system} and  parameters are listed in Table \ref{tab_2_bus_system}.

\begin{table}[h]
\scriptsize
 \begin{center}
   \caption{Parameters for the 2-bus system, where the units for all parameters are p.u. except for $\Delta\theta_{12}^\text{min},\Delta\theta_{12}^\text{max}$ (in degree) and $c_1,c_0$, which are generator cost coefficients.}\label{tab_2_bus_system}
   	\setlength{\tabcolsep}{1 pt}
	\renewcommand\arraystretch{1.2}
   \begin{tabular}{l l | l l}
     \hline
        Symbol      & Value & Symbol      & Value  \\
        \hline
        $Z_{l},Y^{\text{sh}}_{lij}, Y^{\text{sh}}_{lji}, $ & $0.065+j0.62, j0.9, j0.9$ & $P_{g=2}^\text{min},P_{g=2}^\text{max}$ & $0.0,2.0$ \\
        $S_{d=1},S_{d=2}$ & $0.11+j0.4, 0.9+j0.5$ &$Q_{g=2}^\text{min},Q_{g=2}^\text{max}$ & $-1.0,1.0$\\
        $P_{g=1}^\text{min},P_{g=1}^\text{max}$ & $0.0,2.0$ &$S_{l}^\text{max},I_{l}^\text{max}$ & $0.80,0.8551$\\
        $Q_{g=1}^\text{min},Q_{g=1}^\text{max}$ & $-1.0,1.0$ &$U_i^\text{min}, U_i^\text{max}$ & $0.94,1.10$\\
        $\Delta\theta_{ij}^\text{min},\Delta\theta_{ij}^\text{max}$ & $-30^\circ,30^\circ$ &$c_{g,1},c_{g,2}$ & \makecell[l]{g=1: $20,~0$\\g=2: $-10,~0$}\\
        \hline 
   \end{tabular}
 \end{center}
\end{table} 

For the studied case, the optimal objective value of the nonlinear AC OPF problem is 11.5816. 
The \textsc{Canonical BIM} SOC model leads to a gap of 105.89\,\% while the gap for the \textsc{Improved BIM/BFM} formulations with implied total current limit is only  2.92\,\%. As shown in Table \ref{tab_2_bus_system}, the implied bound for the total current is,
\begin{IEEEeqnarray}{C}
    I_{lij}^\text{max}=I_{lji}^\text{max}={S_{l}^\text{max}}/{U_i^\text{min}}={S_{l}^\text{max}}/{U_j^\text{min}}=0.8551.
\end{IEEEeqnarray}
% To investigate how the implied total current limit improves the relaxation gap, the total current are further analysed and compared.
For the \textsc{Canonical BIM}, the optimal total current values at the from and to ends are 4.061 and 4.437, respectively.
In the \textsc{Improved BIM} formulation, the total current at both ends of the line are exactly at the limit, i.e., 0.8551. 
We note that the current limit  really puts a lot of pressure on the feasible space of the power flow losses in the individual branches.

% Explain the following
% \begin{itemize}
% \item Canonical BIM SOC relaxation with only apparent power limits
% \item values for implied total current limit
% \item how this total current constraint is violated in the optimal solution of the canonical BIM SOC relaxation, and therefore changes the gap when we add it
% \end{itemize}

\subsection{PG Lib Results}
Table \ref{tab_case_results} lists the results of numerical experiments.
We focus on a selection of networks where the proposed improvements close the gap significantly, i.e. more than 0.02\,\%.

\begin{table}[tbh]
\setlength{\tabcolsep}{2pt} % Default value: 6pt
 \begin{center}
   \caption{PGLib\_OPF\_ case results. 
   Parentheses indicate \textsc{Ipopt} output status ``\texttt{ALMOST\_LOCALLY\_SOLVED}'' instead of ``\texttt{LOCALLY\_SOLVED}''.}\label{tab_case_results}
   \begin{tabular}{l r r r r r l l l l}
     \hline
Case      & AC-polar \cite{Babaeinejadsarookolaee2019} & \multicolumn{2}{c}{\textsc{BIM can.} } & \multicolumn{2}{c}{\textsc{BIM impr.} }  & \multicolumn{2}{c}{\textsc{BFM impr.} }  \\ 
          & obj ($\$/h$)  & t (s)  & gap (\%)  & t (s) & gap (\%)  &  t (s) & gap (\%)   \\ 
\hline
39\_epri  & 	138415.56  &  	$< 1$  &  	0.55 &  	$< 1$  &  	0.35 &  	$< 1$  &  	0.35 \\
162\_ieee\_dtc &  	108075.64  &  	1 &  	5.94 &  	1 &  	4.73 &  	1 &  	4.73 \\
588\_sdet &  	313139.78  &  	3 &  	2.14 &  	4 &  	1.91 &  	3 &  	1.91 \\
1888\_rte &  	1402530.82 &  	201 &  	2.04 &  	198	 &  0.97 &  	23 &  	0.97 \\
4661\_sdet &  	2251344.07	  &  	48	 &  1.98 &  	65	 &  1.89 &  	48 &  	1.89 \\
6468\_rte &  	2069730.15	 &  	119 &  	1.12 &  	(207) &  	(0.98) &  	273 &  	0.98 \\
\hline
   \end{tabular}
 \end{center}
\end{table} 

The authors propose the following interpretation of why the total current limits can tighten the gap, despite not using new information (as it would be to use independent information on ampacity). 
In the optimal solution of AC OPF, the implied current limits can only be binding when the apparent power limits are binding as well.
However due to the relaxation, the bounds on the total current variables can nevertheless be binding \emph{before} the apparent power bounds are binding.  

Focusing now on the broader set of \textsc{PG Lib} test cases\footnote{Including the cases where there was no significant change in gap.}, the \textsc{Improved BIM} OPF is slower than the \textsc{Canonical BIM}; whereas \textsc{Improved BFM} and  \textsc{Canonical BIM} are very similar in speed. 
Finally, the \textsc{Improved BIM} and \textsc{Improved BFM} have very similar calculation times (on average BFM 1.7\,\% faster).

We note that in the implementation, we did not introduce auxiliary variables for the linking equalities \S\ref{sec_linking_equalities} but instead performed the substitutions. This obviously makes the BIM denser than the BFM from the perspective of the total current limits. Further exploration on whether elimination of these variables helps or hurts performance is needed. 
However,  basing the distinction between BIM and BFM  on the variable space ($S^{\text{tot}}_{lij}, W_i, W_{ij}$ vs $S^{\text{tot}}_{lij}, W_i, L_l$) may then become ambiguous.

\section{Conclusions}
We propose two novel constraint sets, one for implied total current limits, and one for parallel lines, to strengthen the canonical second-order conic BIM and BFM relaxations for networks with $\Pi$-sections as building blocks. 
Two test cases are developed that illustrate the difference w.r.t. SOC BIM and BFM relaxations without these novel constraints.
We show that, by exploiting implied total current limits (using only apparent power limits and voltage bounds), we can tighten the canonical BIM and BFM relaxations. 
Due to the relaxation step, it is possible that the implied total current bounds become binding before the original apparent power bounds. 
% The novel constraints do not seem to have a negative effect on the calculation speed and reliability, so we believe they should be considered broadly. 
The total current limits seem particularly effective when the objective incentivizes network losses, e.g. when generation cost is negative, or dispatching load can solve congestions. 

Finally, we note that OPF with current limits (as opposed to apparent power limits) is an under-explored topic, despite thermal loading being a function of power loss, which is a quadratic function of \emph{current} (magnitude) - not power. 
The proposed loss bounds can  also be used to strengthen network flow relaxations\,\cite{Coffrin2015d}.

% \begin{itemize}
%     \item For parallel lines, the BIM is superior, for lines with nonnegligible losses, BFM is superior. 
%     \item The loss bounds \eqref{eq_loss_bound_1}-\eqref{eq_loss_bound_8} also work in the network flow relaxation \cite{Coffrin2015d}.
%     \item OPF with current limits is under-explored, despite thermal limits being define in terms of power loss, which is a quadratic function of current. It is easy to see that apparent power limits become useless when voltage gets closer to zero. Improved data sets are needed.
%     \item Argue in favor of understanding PG Lib test cases better, and making sure their properties are varied across all things that influence scalability and reliability
% \end{itemize}

\section*{Acknowledgement}
Special thanks to Dr. Hassan Hijazi for alerting us to a typo in the case2\_parallel data. 

\bibliographystyle{IEEEtran}

\begin{thebibliography}{10}
\providecommand{\url}[1]{#1}
\csname url@samestyle\endcsname
\providecommand{\newblock}{\relax}
\providecommand{\bibinfo}[2]{#2}
\providecommand{\BIBentrySTDinterwordspacing}{\spaceskip=0pt\relax}
\providecommand{\BIBentryALTinterwordstretchfactor}{4}
\providecommand{\BIBentryALTinterwordspacing}{\spaceskip=\fontdimen2\font plus
\BIBentryALTinterwordstretchfactor\fontdimen3\font minus
  \fontdimen4\font\relax}
\providecommand{\BIBforeignlanguage}[2]{{%
\expandafter\ifx\csname l@#1\endcsname\relax
\typeout{** WARNING: IEEEtran.bst: No hyphenation pattern has been}%
\typeout{** loaded for the language `#1'. Using the pattern for}%
\typeout{** the default language instead.}%
\else
\language=\csname l@#1\endcsname
\fi
#2}}
\providecommand{\BIBdecl}{\relax}
\BIBdecl

\bibitem{Xu2021}
\BIBentryALTinterwordspacing
S.~Xu, R.~Ma, D.~K. Molzahn, H.~Hijazi, and C.~Josz, ``{Verifying global
  optimality of candidate solutions to polynomial optimization problems using a
  determinant relaxation hierarchy},'' 2021. [Online]. Available:
  \url{http://arxiv.org/abs/2101.00621}
\BIBentrySTDinterwordspacing

\bibitem{Gopinath2020}
S.~Gopinath, H.~L. Hijazi, T.~Weisser, H.~Nagarajan, M.~Yetkin, K.~Sundar, and
  R.~W. Bent, ``{Proving global optimality of ACOPF solutions},''
  \emph{Electric Power Syst. Res.}, vol. 189, no. October 2019, p. 106688,
  2020.

\bibitem{Sundar2018a}
\BIBentryALTinterwordspacing
K.~Sundar, H.~Nagarajan, S.~Misra, M.~Lu, C.~Coffrin, and R.~Bent,
  ``{Optimization-Based Bound Tightening using a Strengthened QC-Relaxation of
  the Optimal Power Flow Problem},'' 2018. [Online]. Available:
  \url{http://arxiv.org/abs/1809.04565}
\BIBentrySTDinterwordspacing

\bibitem{Low2014}
S.~H. Low, ``{Convex relaxation of optimal power flow - part I: formulations
  and equivalence},'' \emph{IEEE Trans. Control Netw. Syst.}, vol.~1, no.~1,
  pp. 15--27, mar 2014.

\bibitem{Babaeinejadsarookolaee2019}
S.~Babaeinejadsarookolaee, A.~Birchfield, R.~D. Christie, C.~Coffrin,
  C.~DeMarco, R.~Diao, M.~Ferris, S.~Fliscounakis, S.~Greene, R.~Huang,
  C.~Josz, R.~Korab, B.~Lesieutre, J.~Maeght, D.~K. Molzahn, T.~J. Overbye,
  P.~Panciatici, B.~Park, J.~Snodgrass, and R.~Zimmerman, ``{The power grid
  library for benchmarking AC optimal power flow algorithms},''
  \emph{[math.OC]}, pp. 1--17, 2019.

\bibitem{coffrin2015distflow}
C.~Coffrin, H.~L. Hijazi, and P.~{Van Hentenryck}, ``{DistFlow extensions for
  AC transmission systems},'' \emph{[Math.OC]}, pp. 1--19, 2015.

\bibitem{Zimmerman2011}
R.~D. Zimmerman, C.~E. Murillo-S{\'{a}}nchez, and R.~J. Thomas, ``{MATPOWER:
  steady-state operations, systems research and education},'' \emph{IEEE Trans.
  Power Syst.}, vol.~26, no.~1, pp. 12--19, 2011.

\bibitem{Coffrin2017}
C.~Coffrin, R.~Bent, K.~Sundar, Y.~Ng, and M.~Lubin, ``{PowerModels.jl: an
  open-source framework for exploring power flow formulations},'' in
  \emph{Power Syst. Comp. Conf.}, vol.~20, Dublin, Ireland, 2018, p.~8.

\bibitem{Coffrin2015c}
C.~Coffrin, H.~Hijazi, and P.~{Van Hentenryck}, ``{Strengthening the SDP
  relaxation of ac power flows with convex envelopes, bound tightening, and
  valid inequalities},'' \emph{IEEE Trans. Power Syst.}, vol.~32, no.~5, pp.
  3549--3558, 2017.

\bibitem{Wachter2006}
A.~W{\"{a}}chter and L.~T. Biegler, ``{On the implementation of primal-dual
  interior point filter line search algorithm for large-scale nonlinear
  programming},'' \emph{Math. Prog.}, vol. 106, no.~1, pp. 25--57, 2006.

\bibitem{Coffrin2015d}
C.~Coffrin, H.~L. Hijazi, and P.~{Van Hentenryck}, ``{Network flow and copper
  plate relaxations for AC transmission systems},'' in \emph{Power Syst. Comp.
  Conf.}, Genoa, 2016, pp. 1--8.

\end{thebibliography}

\end{document}